\documentclass[1p,12pt,number,sort&compress,times]{elsarticle}
\usepackage[utf8]{inputenc}
\usepackage{amsmath,amssymb,amsfonts,amsthm,bbm,bm}
\usepackage{enumitem}
\usepackage{wrapfig}
\usepackage{url}
\usepackage{todonotes}
\setuptodonotes{color=red!40,fancyline}
\usepackage{hyperref}

\newcommand{\RR}{\mathbbm{R}}
\newcommand{\uv}{\mathbf{u}}

\newtheorem{thm}{Theorem}
\newdefinition{rmk}[thm]{Remark}

\begin{document}

\begin{frontmatter}

\title{On the Numerical Integration of Singular Initial and Boundary Value
  Problems for Generalised Lane--Emden and Thomas--Fermi Equations} 
\author[kassel]{Werner M. Seiler}
\ead[url]{http://www.mathematik.uni-kassel.de/~seiler}
\ead{seiler@mathematik.uni-kassel.de}
\author[kassel]{Matthias Sei\ss}
\ead{mseiss@mathematik.uni-kassel.de}
\address[kassel]{Institut f\"ur Mathematik, Universit\"at Kassel, 34132 Kassel,
  Germany}

\begin{abstract}
  We propose a geometric approach for the numerical integration of singular
  initial and boundary value problems for (systems of) quasi-linear differential
  equations. It transforms the original problem into the problem of
  computing the unstable manifold at a stationary point of an associated
  vector field and thus into one which can be solved in an efficient and
  robust manner. Using the shooting method, our approach also works well
  for boundary value problems. As examples, we treat some (generalised)
  Lane--Emden equations and the Thomas--Fermi equation.
\end{abstract}


\end{frontmatter}

\section{Introduction}

The \emph{Lane--Emden equation} was originally derived in astrophysics
\cite[p.~40]{re:gas} and represents a dimensionless form of Poisson's
equation for the gravitational potential of a Newtonian self-gravitating,
spherically symmetric, polytropic fluid (see
\cite{htd:diffint,hkt:stellar,gph:poly} and references therein for a more
detailed discussion):
\begin{equation}\label{eq:lesc}
  u''+\frac{2}{x}u'=-u^{n}
\end{equation}
with $n$ the polytropic index.  Astrophysicists need the solution of the
\emph{initial value problem} $u(0)=1$ and $u'(0)=0$.  Eq.~\eqref{eq:lesc}
is prototypical for ordinary differential equations arising in the
construction of radially symmetric steady state solutions of
reaction-diffusion equations, as the left hand side of \eqref{eq:lesc}
represents the Laplace operator in spherical coordinates.  In an
$N$-dimensional space, the numerator $2$ has to be replaced by $N-1$.  This
leads to \emph{generalised Lane--Emden equations} of the form
\begin{equation}\label{eq:lesg}
  u''+\frac{N-1}{x}u'=h(x,u)\;,
\end{equation}
where the function $h$ represents the reaction term.  Besides the classical
form from astrophysics, we will later consider examples arising in chemical
engineering (biocatalysts) and in physiology (oxygen uptake of cells).
There, one needs the solution of \emph{boundary value problems} with
$u'(0)=0$ and $\alpha u(1)+\beta u'(1)=\gamma$.

Thomas \cite{lht:atom} and Fermi \cite{ef:stat} derived independently of
each other in a statistical model of atoms treating electrons as a gas of
particles a Lane--Emden equation \eqref{eq:lesc} with polytropic index
$n=3/2$ for the electrostatic potential $V(x)$, however with the ``initial
condition'' that $V(x)$ behaves like $1/x$ for $x\to0$.  Writing
$V(x)=u(x)/x$, one obtains the \emph{Thomas--Fermi equation}
\begin{equation}\label{eq:tf}
  u''=\sqrt{u^{3}/x}
\end{equation}
together with the initial condition $u(0)=1$ (see
\cite{dge:fmsma,nhm:tfqm,imt:iap} for more physical and historical details
and \cite{eh:tf1,eh:tf2} for a mathematical analysis).  In addition, one
imposes one of the following three types of boundary conditions:
\begin{subequations}\label{eq:tfbc}
  \begin{align}
    &bu'(b)-u(b)=0\;,\label{eq:tfbccryst}\\
    &\!\!\lim_{x\to\infty}u(x)=0\;,\label{eq:tfbcinf}\\
    &u(a)=0\label{eq:tfbcion}
  \end{align}
\end{subequations}
with $0<a,b<\infty$ given positions.  The infinite case \eqref{eq:tfbcinf}
occurs only for a critical value $\omega\approx-1.588\dots$ of the initial
slope $u'(0)$ and represents physically an isolated neutral atom.  For
larger initial slopes, one can prescribe the boundary condition
\eqref{eq:tfbccryst} and obtains solutions going through a minimum and then
growing rapidly.  Physically, such solutions are relevant for certain
crystals.  The boundary condition \eqref{eq:tfbcion} leads to solutions
with a smaller initial slope and represent physically ions with radius $a$.

Numerical methods from textbooks cannot be directly applied here, as all
considered equations are singular at $x=0$ and at least one
initial/boundary condition is imposed there.  In the vast literature on the
numerical integration of Lane--Emden and Thomas--Fermi equations, three
different types of approaches prevail.  For initial value problems,
astrophysicists apply a very simple approach based on a series expansion of
the solution to get away from the singularity and then some standard
integrator \cite[Sect.~7.7.2]{hkt:stellar} (see also \cite{gph:7dle}).  For
boundary value problems, collocation methods are popular, as they are
easily adapted to the singularity \cite{rs:numsbvp}.  Finally, various
kinds of semi-analytic expansions like Adomian decomposition have been
adapted to the singularity (see the references given below and references
therein).

We propose here a new and rather different alternative.  In the geometric
theory of differential equations \cite{via:geoode,aor:poincare}, one
associates with any implicit ordinary differential equation a vector field
on a higher-dimensional space such that the graphs of prolonged solutions
of the implicit equation are integral curves of this vector field.  Most of
the literature on singularity theory is concerned with fully implicit
equations.  However, in real life applications quasi-linear equations like
the Lane--Emden equations prevail.  In \cite{wms:singbif,wms:quasilin}, the
authors showed that such equations possess a special geometry allowing us
to work in a lower order.  Singularities, now called impasse points, are
typically stationary points of the associated vector field.  If there is a
unique solution, its prolonged solution graph is the one-dimensional
unstable manifold of this stationary point.  Such an unstable manifold can
numerically be computed very robustly.  In \cite{bss:numvis}, we already
sketched this possibility to exploit ideas from singularity theory for
numerical analysis.  Here, we want to demonstrate for concrete problems of
practical relevance that it is easy to apply and efficiently provides
accurate results.

The paper is structured as follows.  In the next section, we recall the
necessary elements of the geometric theory of differential equations and
how one can translate an implicit problem into an explicit one.
Section~\ref{sec:le} is then devoted to the application of these ideas to
(generalised) Lane--Emden equations and to the numerical solution of some
concrete problems from the literature.  In Section~\ref{sec:tf} we discuss
the Thomas--Fermi equation by first reducing it via a transformation
introduced by Majorana and then applying the geometric theory.  We compare
the obtained numerical results with some high precision calculations from
the literature.  Finally, some conclusions are given.

\section{Geometric Theory of Ordinary Differential Equations} 
\label{sec:geotheo}

We use a differential geometric approach to differential equations.  It is
beyond the scope of this article to provide deeper explanations of it; for
this we refer to \cite{wms:invol} and references therein.  For notational
simplicity, we concentrate on the scalar case; the extension to systems
will be briefly discussed at the end.  Similarly, we restrict here to
second-order equations, but equations of arbitrary order can be treated in
an analogous manner.

We consider a fully implicit differential equation of the form
\begin{equation} \label{eq:impl2}
  F(x,u,u',u'')=0\;.
\end{equation}
In the second-order\emph{ jet bundle} $\mathcal{J}_{2}$ (intuitively
expressed, this is simply a four-dimensional affine space with coordinates
called $x,u,u',u''$), this equation defines a hypersurface
$\mathcal{R}_{2}\subset\mathcal{J}_{2}$ which represents our geometric
model of the differential equation.  We will assume throughout that
$\mathcal{R}_{2}$ is actually a submanifold.

Given a function $\psi(x)$, we may consider its \emph{graph} as a curve in
the jet bundle $\mathcal{J}_{0}$ of order zero, i.\,e.\ the $x$-$u$ space,
given by the map $x\mapsto\bigl(x,\psi(x))$.  Assuming that $\psi$ is at
least twice differentiable, we can \emph{prolong} this curve to a curve in
$\mathcal{J}_{2}$ defined by the map
$x\mapsto \bigl(x,\psi(x),\psi'(x),\psi''(x)\bigr)$.  The function $\psi$
is a solution of \eqref{eq:impl2}, if and only if this curve lies
completely in the hypersurface $\mathcal{R}_{2}$.

In an \emph{initial value problem} for the implicit equation
\eqref{eq:impl2}, one prescribes a point
$\rho=(y,u_{0},u_{1},u_{2})\in\mathcal{R}_{2}$ and asks for solutions such
that $\rho$ lies on their prolonged graphs.  Note that opposed to explicit
problems, we must also specify the value~$u_{2}$, as the algebraic equation
$F(y,u_{0},u_{1},u'')=0$ may have several (possibly infinitely many)
solutions and thus may not uniquely determine $u_{2}$.

A key ingredient of the geometry of jet bundles is the contact structure.
In the case of $\mathcal{J}_{2}$, the \emph{contact distribution}
$\mathcal{C}^{(2)}$ is spanned by the two vector fields
\begin{equation}\label{eq:condist2}
  C_{\mathrm{trans}}=\partial_{x}+u'\partial_{u}+u''\partial_{u'}\;,\quad
  C_{\mathrm{vert}}=\partial_{u''}\;.
\end{equation}
A curve $x\mapsto \bigl(x,\psi_{0}(x),\psi_{1}(x),\psi_{2}(x)\bigr)$ in
$\mathcal{J}_{2}$ is a prolonged graph (i.\,e.\ $\psi_{1}=\psi_{0}'$ and
$\psi_{2}=\psi_{0}''$), if and only if all its tangent vectors lie in the
contact distribution.

The \emph{Vessiot distribution} $\mathcal{V}[\mathcal{R}_{2}]$ of
\eqref{eq:impl2} is that part of the tangent space of $\mathcal{R}_{2}$
which also lies in the contact distribution $\mathcal{C}^{(2)}$.  Writing
$X=a C_{\mathrm{trans}}+ b C_{\mathrm{vert}}$ for a general vector in the
contact distribution, $X$ lies in the Vessiot distribution, if and only if
its coefficients $a,b$ satisfy the linear equation
\begin{equation}\label{eq:vdimpl2}
  \bigl(F_{x}+u'F_{u}+u''F_{u'}\bigr)a+F_{u''}b=0\;.
\end{equation}
A \emph{singularity} is a point
$\rho=(y,u_{0},u_{1},u_{2})\in\mathcal{R}_{2}$ such that $F_{u''}(\rho)=0$.
One speaks of a \emph{regular singularity}, if the coefficient of $a$ in
\eqref{eq:vdimpl2} does not vanish at $\rho$, and of an \emph{irregular
  singularity}, if it does.  Outside of irregular singularities, the
Vessiot distribution is one-dimensional and locally spanned by the vector
field
\begin{equation}\label{eq:Xvd}
  X=F_{u''}C_{\mathrm{trans}}-
  \bigl(F_{x}+u'F_{u}+u''F_{u'}\bigr)C_{\mathrm{vert}}
\end{equation}
(note that $X$ is defined only on the submanifold
$\mathcal{R}_{2}\subset\mathcal{J}_{2}$).  The prolonged graph of any
solution of \eqref{eq:impl2} must be integral curves of this vector field.
The converse is not necessarily true in the presence of singularities.

At regular singularities, the vector field $X$ becomes vertical.
Generically, only one-sided solutions exist at such points and if two-sided
solutions exist, then their third derivative will blow up
\cite[Thm.~4.1]{wms:aims}.  At irregular singularities, typically several
(possibly infinitely many) solutions exist.  In \cite{lrss:gsade} it is
shown how for arbitrary systems of ordinary or partial differential
equations with polynomial nonlinearities all singularities can be
automatically detected.

Irregular singularities are stationary points of $X$.  Prolonged solution
graphs through them are one-dimensional invariant manifolds.  Any
one-dimensional (un)stable or centre manifold (with transversal tangent
vectors) at such a stationary point defines a solution.  For
higher-dimensional invariant manifolds, one must study the induced dynamics
on them to identify solutions.  In any case, we note that the numerical
determination of invariant manifolds at stationary points is a well-studied
topic -- see e.\,g.\ \cite{bk:invman,ep:teim}.

In general, the direct numerical integration of \eqref{eq:impl2} faces some
problems, if it is not possible to solve (uniquely) for $u''$, and
typically breaks down, if one gets too close to a singularity.  The
geometric theory offers here as alternative the numerical integration of
the dynamical system defined by the vector field~$X$.  Thus an implicit
problem is transformed into an explicit one!  The price one has to pay is
an increase of the dimension: while \eqref{eq:impl2} is a scalar equation
(but second-order), the vector field $X$ lives on the three-dimensional
manifold $\mathcal{R}_{2}$ in the four-dimensional jet bundle
$\mathcal{J}_{2}$ (more generally, a scalar equation of order $q$ leads to
a vector field on a $(q-1)$-dimensional manifold).

The key difference is, however, that we obtain a parametric solution
representation.  We work now with the explicit autonomous
system\footnote{Strictly speaking, we are dealing here with a
  three-dimensional system, as $X$ lives on the three-dimensional
  manifold~$\mathcal{R}_{2}$.  As we do not know a parametrisation of
  $\mathcal{R}_{2}$, we must work with all four coordinates of
  $\mathcal{J}_{2}$.  One could augment \eqref{eq:sysvd} by its first
  integral $F(x,u,u',u'')=0$ and enforce it during a numerical integration,
  but in our experience this is not necessary.}
\begin{equation}\label{eq:sysvd}
  \begin{aligned}
    \frac{dx}{ds} &= F_{u''}\;, &
    \frac{du}{ds} &= u'F_{u''}\;,\\
    \frac{du'}{ds} &= u''F_{u''}\;,\qquad &
    \frac{du''}{ds} &= -F_{x} - u'F_{u} - u''F_{u'}\;,
  \end{aligned}
\end{equation}
where $s$ is some auxiliary variable used to parametrise the integral
curves of $X$.  A solution of it will thus be a curve
$s\mapsto\bigl(x(s),u(s),u'(s),u''(s)\bigr)$ on
$\mathcal{R}_{2}\subset\mathcal{J}_{2}$.  A numerical integration will
provide a discrete approximation of this curve. 

In applications, \emph{quasi-linear} equations prevail.  We restrict here
even to semi-linear differential equations of the form
\begin{equation}\label{eq:sql2e}
  F(x,u,u',u'')=g(x)u''-f(x,u,u')=0\;,
\end{equation}
with smooth functions $f$, $g$, as both the Lane--Emden and the
Thomas--Fermi equation can be brought into this form.  A point
$(y,u_{0},u_{1},u_{2})\in\mathcal{R}_{2}$ is then a singularity, if and
only if $g(y)=0$.

As first shown in \cite{wms:singbif} and later discussed in more details in
\cite{wms:quasilin}, quasi-linear equations possess their own special
geometry, as it is possible to project the Vessiot distribution to the jet
bundle of one order less, i.\,e.\ in our case to the first-order jet bundle
$\mathcal{J}_{1}$ with coordinates $(x,u,u')$.  Projecting the vector field
$X$ defined by \eqref{eq:Xvd} with $F$ as in \eqref{eq:sql2e} to
$\mathcal{J}_{1}$ yields the vector field
\begin{equation}\label{eq:Yvd}
  Y=g(x)\partial_{x}+g(x)u'\partial_{u}+f(x,u,u')\partial_{u'}\;.
\end{equation}
It is only defined on the canonical projection of $\mathcal{R}_{2}$ to
$\mathcal{J}_{1}$ which may be a proper subset.  Assuming that $f$, $g$ are
defined everywhere on $\mathcal{J}_{1}$, we analytically extend $Y$ to all
of $\mathcal{J}_{1}$ and replace \eqref{eq:sysvd} by the three-dimensional
system
\begin{equation}\label{eq:sysvdql}
  \frac{dx}{ds} = g(x)\;,\quad
  \frac{du}{ds} = g(x)u'\;,\quad
  \frac{du'}{ds} = f(x,u,u')\;.
\end{equation}
The first equation is decoupled and can be interpreted as describing a
change of the independent variable, but we will not pursue this point of
view.

A point $\rho=(y,u_{0},u_{1})\in\mathcal{J}_{1}$ is an \emph{impasse point}
for \eqref{eq:sql2e}, if the vector field $Y$ is not transversal at $\rho$,
i.\,e.\ if its $x$-component vanishes.  Here, this is equivalent to
$g(y)=0$.  We call $\rho$ a \emph{proper} impasse point, if
$\mathcal{R}_{2}$ contains points which project on $\rho$; otherwise,
$\rho$ is \emph{improper}.  Here, proper impasse points are obviously
stationary points of $Y$ or \eqref{eq:sysvdql}, respectively.  Prolonged
graphs of solutions of \eqref{eq:sql2e} are one-dimensional invariant
manifolds of $Y$ (or \eqref{eq:sysvdql}, resp.) and again the converse is
not necessarily true.  In \cite{wms:quasilin}, we proved geometrically the
following result (a classical analytic proof for the special case $g(x)=x$
can be found in \cite{jfl:singivp}).

\begin{thm}\label{thm:exunisql2}
  Consider \eqref{eq:sql2e} for $f,g$ smooth together with the initial
  conditions $u(y)=u_{0}$ and $u'(y)=u_{1}$ where $g(y)=0$ and
  $f(y,u_{0},u_{1})=0$.  If $\delta=g'(y)$ and
  $\gamma=f_{u'}(y,u_{0},u_{1})$ are both non zero and of opposite sign,
  then the initial value problem possesses a unique smooth solution.
\end{thm}

Under the made assumptions, the initial point $\rho=(y,u_{0},u_{1})$ is a
proper impasse point of \eqref{eq:sql2e}.  One readily verifies that the
Jacobian~$J$ of $Y$ at $\rho$ has the eigenvalues $\delta$, $0$ and
$\gamma$ and thus we find three one-dimensional invariant manifolds at
$\rho$: the stable, the unstable and the centre manifold.\footnote{The
  centre manifold is here unique, as there exists a whole curve of
  stationary points \cite{js:cm}.}  Without loss of generality, we assume
that $\delta>0$ (otherwise we multiply \eqref{eq:sql2e} by $-1$).  It is
then shown in \cite{wms:quasilin} that the prolonged graph of the unique
solution is the unstable manifold and thus at $\rho$ it is tangent to the
eigenvector of $J$ for $\delta$.

\begin{rmk}
  The extension to implicit systems $\mathbf{F}(x,\uv,\uv',\uv'')=0$ is
  straightforward.  Assuming that the unknown function $\uv$ is vector
  valued, $\uv\colon\mathcal{I}\subseteq\RR\to\RR^{n}$, the jet
  bundle~$\mathcal{J}_{2}$ is $(3n+1)$-dimensional and the contact
  distribution $\mathcal{C}^{(2)}$ is generated by the $n+1$ vector fields
  $C_{\mathrm{trans}}=\partial_{x}+\uv'\cdot\partial_{\uv}+
  \uv''\cdot\partial_{\uv'}$ and
  $\mathbf{C}_{\mathrm{vert}}=\partial_{\uv''}$, where the dot denotes the
  standard scalar product.  Again the Vessiot distribution is generically
  one-dimensional and the coefficients of a vector field $X$ spanning it
  are readily determined by solving a linear system of equations.
  Numerical integration of $X$ allows us to approximate solutions of the
  given system.

  We restrict to semi-linear first-order systems of the form
  $g(x)\uv'=\mathbf{f}(x,\uv)$ with $g$ still a scalar function.  For
  initial conditions $\uv(y)=\uv_{0}$ with $g(y)=0$ and
  $\mathbf{f}(y,\uv_{0})=0$, we introduce $\delta=g'(y)$ (assuming
  $\delta>0$) and the Jacobian $\Gamma=\mathbf{f}_{\uv}(y,\uv_{0})$.  In
  \cite{wms:quasilin2}, it is shown that if all eigenvalues of $\Gamma$
  have a negative real part, then the initial value problem has a unique
  smooth solution.  A classical analytical proof was given by Vainikko by
  first studying extensively the linear case \cite{gv:linsingode} and then
  extending to the nonlinear one \cite{gv:nonlinsingode}.  In the geometric
  approach, one sees again that the graph of the solution is a
  one-dimensional unstable manifold of the vector field $Y$ spanning the
  projected Vessiot distribution.
\end{rmk}

\section{(Generalised) Lane--Emden Equations}
\label{sec:le}

\subsection{Geometric Treatment}
\label{sec:legeo}

If we consider the generalised Lane--Emden equation \eqref{eq:lesg}, then
one obtains after multiplication by $x$ the special case of
\eqref{eq:sql2e} given by
\begin{equation}\label{eq:gfle}
  g(x)=x\;,\qquad f(x,u,u')=xh(x,u)-(N-1)u'\;,
\end{equation}
where we always assume $N>1$.  For arbitrary initial conditions
$u(0)=u_{0}$ and $u'(0)=u_{1}$, we find that $\delta=1$ and $\gamma=-(N-1)$
are nonzero and of opposite sign.  The initial point $\rho=(0,u_{0},u_{1})$
is a \emph{proper} impasse point, if and only if $u_{1}=0$.  In this case,
Theorem~\ref{thm:exunisql2} asserts the existence of a unique smooth
solution.

The projected Vessiot distribution is spanned by the vector field
\begin{equation}\label{eq:sleY}
  Y=x\partial_{x}+xu'\partial_{u}+\bigl[xh(x,u)-(N-1)u'\bigr]\partial_{u'}\;.
\end{equation}
For $u_{1}\neq0$, no solution can exist.  Indeed, the vector field $Y$ has
then no stationary point and the unique trajectory through the initial
point $\rho=(0,u_{0},u_{1})$ is the vertical line
$s\mapsto(0,u_{0},u_{1}+s)$ which does not define a prolonged graph.

We thus assume $u_{1}=0$, which unsurprisingly is the case in all
applications of \eqref{eq:lesg} in the literature.  Independent of the
value of $u_{0}$, the initial point $\rho=(0,u_{0},0)$ is a stationary
point of the vector field $Y$.  The Jacobian of $Y$ at $\rho$ is
\begin{equation}\label{eq:sleJ}
  J=
  \begin{pmatrix}
    1 & 0 & 0\\
    0 & 0 & 0\\
    h(0,u_{0}) & 0 & -(N-1)
  \end{pmatrix}\;.
\end{equation}
Its eigenvalues are $1$, $0$ and $-(N-1)$.  Relevant for us is only the
eigenvector to the eigenvalue $1$, as it is tangential to the unstable
manifold.  It is given by $\mathbf{v}=\bigl(N,0,h(0,u_{0})\bigr)^{T}$.

For the numerical solution of our given initial value problem, we integrate
the vector field $Y$ for the initial data
$\bigl(x(0),u(0),u'(0)\bigr)^{T}=\bigl(0,u_{0},0\bigr)^{T}+\epsilon\mathbf{v}$
with some small $\epsilon>0$.  The concrete value of $\epsilon$ is not very
relevant.  As the exact solution corresponds to the unstable manifold, any
error is automatically damped by the dynamics of $Y$.  In our experiments,
we typically used $\epsilon=10^{-3}$ or $\epsilon=10^{-4}$.

We can easily extend this approach to coupled systems of the form
\begin{equation}\label{eq:lecs}
  \uv''+\frac{N-1}{x}\uv'=\mathbf{h}(x,\uv)\;,
\end{equation}
where $\uv$ is a vector valued function and the coupling occurs solely
through the reaction terms.  If $\uv$ is a $d$-dimensional vector, then the
dimension of the first-order jet bundle $\mathcal{J}_{1}$ is $2d+1$.  Thus
\eqref{eq:sysvdql} becomes a system of this dimension:
\begin{equation}\label{eq:sysvdlecs}
  \frac{dx}{ds} =x\;,\quad
  \frac{d\uv}{ds}=x\uv'\;,\quad
  \frac{d\uv'}{ds}=x\mathbf{h}(x,\uv)-(N-1)\uv'\;.
\end{equation}
By the same arguments as in the scalar case, we restrict to the initial
condition $\uv'(0)=\mathbf{0}$ so that the initial point
$\rho=(0,\uv_{0},\mathbf{0})$ is again a proper impasse point.  The
Jacobian at $\rho$ is a block form of \eqref{eq:sleJ}:
\begin{equation}\label{eq:lecsJ}
    J=
  \begin{pmatrix}
    1 & \mathbf{0}^{T} & \mathbf{0}^{T} \\
    \mathbf{0} & 0_{d} & 0_{d} \\
    \mathbf{h}(0,\uv_{0}) & 0_{d} & -(N-1)E_{d}
  \end{pmatrix}\;,
\end{equation}
where $0_{d}$ and $E_{d}$ denote the $d\times d$ zero and unit matrix,
resp.  We still have $1$ as a simple eigenvalue, whereas the eigenvalues
$0$ and $-(N-1)$ have both the algebraic multiplicity $d$.  The
$d$-dimensional stable and centre manifolds are again vertical and
irrelevant for a solution theory.  But we still find a one-dimensional
unstable manifold corresponding to the prolonged graph of the unique
solution.  It is tangential to the vector
$\mathbf{v}=\bigl(N,\mathbf{0}^{T},\mathbf{h}(0,\uv_{0})^{T}\bigr)^{T}$ and
as in the scalar case we use as initial data for its determination the
point $\bigl(0,\uv_{0}^{T},\mathbf{0}^{T}\bigr)^{T}+\epsilon \mathbf{v}$.

\subsection{Numerical Results}
\label{sec:lenum}

As our main goal consists of showing how easy the numerical integration of
singular problems becomes with our geometric approach, we did not develop
any sophisticated production code.  We performed all our computations with
the built-in numerical capabilities of \textsc{Maple}.  We used most of the
time the \texttt{dsolve/numeric} command with its standard settings,
i.\,e.\ a Runge--Kutta--Fehlberg pair of order 4/5 is applied with a
tolerance of $10^{-6}$ for the relative error and $10^{-7}$ for the
absolute error.

Our geometric ansatz does not determine approximations
$u_{n}\approx u(x_{n})$ of the solution $u(x)$ on a discrete mesh
$(x_{n})$, but approximations $x_{n}=x(s_{n})$ and $u_{n}=u(s_{n})$ for a
\emph{parametric} representation $\bigl(x(s),u(s)\bigr)$ of the graph of
the solution.  Hence, for computing an approximated solution value
$u(\bar{x})$, one must first determine a parameter value $\bar{s}$ such
that $x(\bar{s})\approx\bar{x}$.  This can easily be accomplished either
with a nonlinear solver or with a numerical integrator with event handling.
We used the latter option in most of our experiments.

For boundary value problems, we applied the shooting method which worked
very well.  Since \textsc{Maple} provides no built-in command for it, we
wrote our own simple version.  In scalar problems, we solved the arising
nonlinear equation most of the time with the Steffensen method (with Aitken
$\Delta^{2}$ acceleration).  As our equations are dimension-free, suitable
starting values were easy to find: typically, $u(x)$ varied between $0$ and
$1$ and we simply chose the midpoint $0.5$.

We encountered difficulties only in the simulation of a biocatalyst.  For
some parameter values, the correct initial value was very close to zero and
the Steffensen iterations produced sometimes intermediate approximations
which were negative and for which the numerical integration became
meaningless.  Here we resorted to a simple bisection method.

For Lane--Emden systems, we used the Newton method for the arising
nonlinear systems.  The Jacobian was determined via the variational
equation of the differential system.  Thus for an $n$-dimensional
differential system where $k<n$ initial conditions have to be determined
via shooting, we had to solve an additional $kn$-dimensional linear
differential system with variable coefficients.

\subsubsection{Scalar Lane--Emden Equations}
\label{sec:slee}

We consider scalar Lane--Emden equations of the generalised form
\begin{subequations}\label{eq:sleall}
\begin{equation}\label{eq:sle}
  u''+\frac{m}{x}u'=f(x,u)
\end{equation}
together with either the initial conditions
\begin{equation}\label{eq:sleic}
  u(0)=u_{0}\;,\quad u'(0)=0 
\end{equation}
or the boundary conditions
\begin{equation}\label{eq:slebc}
  u'(0)=0\;,\quad \alpha u(1)+\beta u'(1)=\gamma\;.
\end{equation}
\end{subequations}
Chawla and Shivakumar \cite{cs:exssbvp} proved for boundary value problems
with $\alpha=1$ and $\beta=0$ an existence and uniqueness theorem under the
following assumption on the right hand side $f(x,u)$: the supremum $M$ of
the negative partial derivative $-f_{u}(x,u)$ on $[0,1]\times\RR$ must be
less than the first positive root $t_{1}$ of the Bessel function
$J_{(m-1)/2}(\sqrt{t})$ (in the frequent case $m=2$, we thus need
$M<\pi^{2}$).

The numerical integration of \eqref{eq:sle} has been studied by many
authors using many different approaches; we refer to \cite{pdrg:lehc} for
an overview of many works before 2010.  We will discuss three different
situations: (i) initial value problems in astrophysics, (ii) Dirichlet
boundary value problem in chemical engineering and (iii) mixed boundary
value problems in physiology.

\paragraph{Initial Value Problems from Astrophysics}
In the classical Lane--Emden equations, one has $m=N-1$ with $N$ the space
dimension and $f(x,u)=-u^{n}$.  The solutions for $u_{0}=1$ are known as
\emph{polytropes}.  Physically meaningful is the range $0\leq n <5$ (with
$n$ not necessarily an integer).  For three polytropic indices, namely
$n=0,1,5$, exact solutions are known \cite[Sect.~2.3]{gph:poly}.  Of
physical relevance are in particular the first zero $\xi_{1}$ of $u$
(corresponding to the scaled radius of the sphere) and the value of
$u'(\xi_{1})$ (e.\,g.\ the ratio of the central density to the mean density
is given by $r=-\xi_{1}/3u'(\xi_{1})$).

\begin{figure}[ht]
  \centering
  \includegraphics[height=6cm]{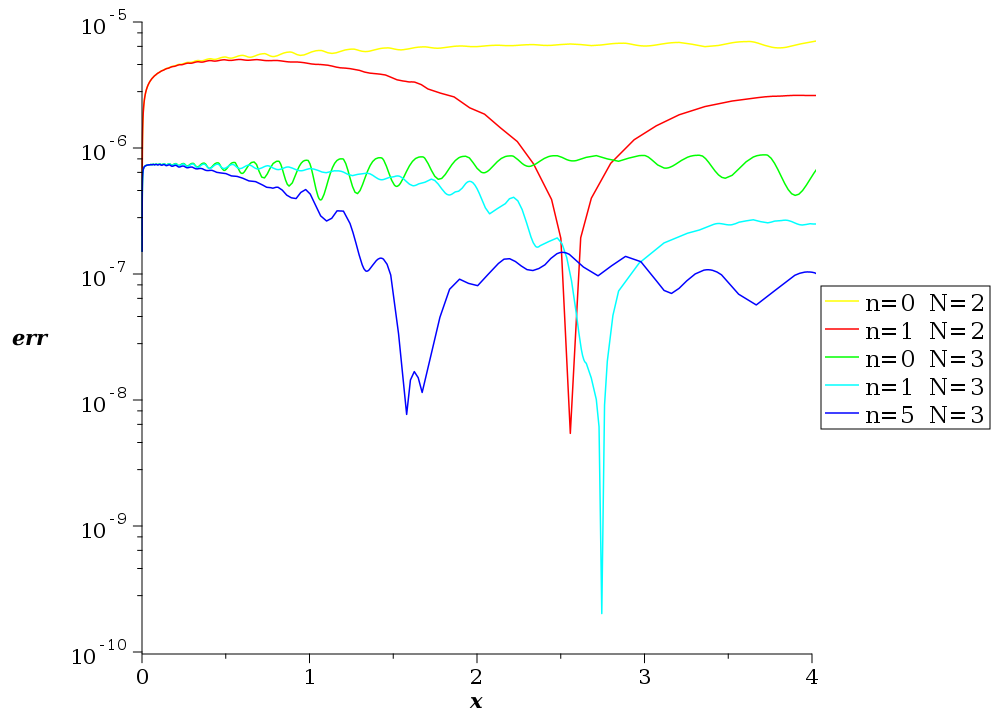}
  \caption{Logarithmic plot of absolute deviation from exact solution for
    some polytropes.}
  \label{fig:polytrope}
\end{figure}

We numerically solved the Lane--Emden equations by integrating the
dynamical system \eqref{eq:sysvdql} with $f,g$ given by \eqref{eq:gfle}.
As concrete test cases, we used some polytropic cylinders and spheres where
the exact solutions are known.  Figure~\ref{fig:polytrope} shows the
observed errors in logarithmic scale.  Obviously, the results are within
the expected range for the default settings of \textsc{Maple}'s numerical
integrator. 

\begin{wraptable}{r}{0.45\linewidth}
  \centering
  \begin{tabular}{rrr}
    \hline
    \multicolumn{1}{c}{$N,\ n$} & \multicolumn{1}{c}{$\xi_{1}$} &
    \multicolumn{1}{c}{$r$} \\
    \hline
    $2,\ 0$ & \ $3.2\cdot 10^{-6}$ & $4.0\cdot 10^{-7\phantom{0}}$\\
    $2,\ 2$ & $4.2\cdot 10^{-7}$ & $5.7\cdot 10^{-6\phantom{0}}$\\
    $3,\ 0$ & $4.3\cdot 10^{-7}$ & $2.2\cdot 10^{-10}$\\
    $3,\ 1$ & $1.5\cdot 10^{-7}$ & $9.3\cdot 10^{-7\phantom{0}}$\\
    \hline
  \end{tabular}
  \caption{Relative errors for first zero $\xi_{1}$ and density ratio $r$
    for the cases with $\xi_{1}<\infty$.}
  \label{tab:polytrope}
\end{wraptable}
Our approach also determines approximations $u_{n}'=u'(s_{n})$ for the
first derivatives of the solution, as the integral curves of the vector
field $Y$ define a parametrisation $\bigl(x(s),u(s),u'(s)\bigr)$ of the
solution and its first derivative.  We use this to approximate also the
quantities $\xi_{1}$ and $r$.  Table~\ref{tab:polytrope} exhibits their
relative errors compared with the exact solution for those cases where
$\xi_{1}$ is finite.  Again, the observed accuracy corresponds well to the
settings of the numerical integrator.

\paragraph{Boundary Value Problems for Biocatalysts}
In chemical engineering, the Lane--Emden equation arises in the analysis of
diffusive transport and chemical reactions of species inside a porous
catalyst pellet \cite[\S6.4]{ra:acr} with boundary conditions of the form
\eqref{eq:slebc} with $\alpha=\gamma=1$ and $\beta=0$.  Flockerzi and
Sundmacher \cite{fs:ledfm} considered the case $m=2$ and
$f(x,u)=\phi^{2}u^{n}$ for a single species obeying Fick's law with
constant diffusivity and power-law kinetics (the constant $\phi^{2}$ is the
Thiele modulus describing the ratio of surface reaction rate to diffusion
rate).  As this corresponds up to a sign exactly to the above considered
polytropes, we omit concrete calculations and only note that
\cite{fs:ledfm} also provides a nice geometric proof of the existence of a
unique solution of this particular boundary value problem which,
unfortunately, seems not be extendable to other functions $f$.

\begin{wrapfigure}{r}{0.45\linewidth}
  \centering
  \includegraphics[width=6cm]{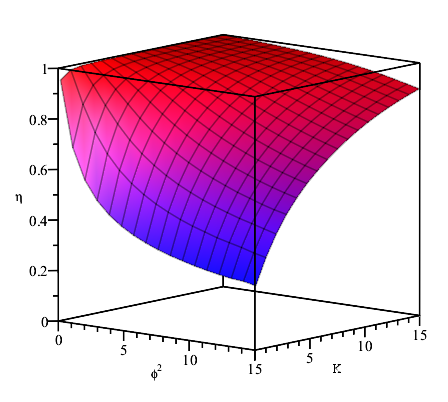}
  \caption{Dependency of the effectiveness factor $\eta$ on Thiele modulus
    $\phi^2$ and dimensionless Michaelis--Menten constant
    $K$.}\label{fig:biocat}
\end{wrapfigure}

Using a Michaelis--Menten kinetics for a biocatalyst, one obtains right
hand sides like $f(x,u)=9\phi^{2}\tfrac{u}{1+K u}$, where $\phi$ is again
the Thiele modulus and $K$ the dimensionless Michaelis--Menten constant
(see \cite{pvr:irk} for some further variants).  This model was analysed by
a homotopy perturbation method in \cite{ass:ibp}.  A quantity relevant for
engineers is the effectiveness factor which is here given by
$\eta=\tfrac{K+1}{3\phi^{2}}u'(1)$.  A numerical study of the dependency of
$\eta$ on $\phi^{2}$ and $K$ leads to the surface shown in
Fig.~\ref{fig:biocat} based on a $17\times17$ grid, i.\,e.\ on the
numerical solution of 289 boundary value problems with different
combinations of parameter values.  As indicated above, we had to use here a
bisection method for locating the right initial value.  Bisecting until an
interval length of $10^{-5}$ was reached, the whole computation required
only 2--3sec on a laptop (equipped with eight Intel Core i7-11370H (11th
generation) working with 3.3GHz and 16GB of RAM running \textsc{Maple 2022}
under Windows 11).

\textsc{Matlab}'s solvers \texttt{bvp4c} and \texttt{bvp5c} are finite
difference methods based on a three- and four-stage, resp., Lobatto IIIa
collocation formulae and provide a special option for the type of
singularity appearing in Lane--Emden equations \cite{sk:bvp4c,sk:bvp5c}.
However, it turned out to be nontrivial to determine a plot like
Fig.~\ref{fig:biocat} with them, as for some parameter values they react
rather sensitive to the required initial guess.  Using a simple constant
function lead sometimes either to completely wrong solutions or the
collocation equations could not be solved.  We then computed one solution
with ``harmless'' parameter values and used it as initial guess for all
other parameter values.  But the computations required with 5--6sec about
twice as much time as our approach.

An alternative approach consists in transforming the problem into a
reaction-diffusion equation by adding a time derivative.  The desired
solution of our boundary value problem arises then as asymptotic for long
times.  \textsc{Matlab} provides here with \texttt{pdepe} a specialised
solver admitting again our type of singularity.  It employs a method for
parabolic partial differential equations proposed by Skeel and Brezins
\cite{sb:sdpe} using a spatial discretisation derived with a Galerkin
approach. Here, one does not need an initial guess and it turns out that a
steady state is reached very rapidly (already $t=1$ is sufficient).  But
one needs an additional interpolation with \texttt{pdeval} to determine
derivative values.  Furthermore, the computation time for a plot like
Fig.~\ref{fig:biocat} increases significantly to about 17sec.\footnote{This
  approach was also used by the authors of \cite{ass:ibp} to compute
  reference solutions.  However, the plots presented there do not agree
  with our results. As they provided a listing of their \textsc{Matlab}
  code, we could repeat their numerical experiments and obtained the same
  results as with our method and not what they show in their paper.}

\paragraph{Mixed Boundary Conditions for a Physiological Model}
The same differential equation is used to model the steady state oxygen
diffusion in a spherical cell with Michaelis-Menten uptake kinetics
\cite{shl:oxygen,me:oxygen}, $m=2$ and $f(x,u)=\tfrac{a u}{u+K}$, but with
mixed boundary conditions \eqref{eq:slebc} where $\alpha=\gamma$,
$\beta=1$.  Hiltmann and Lory \cite{hl:oxygen} proved explicitly the
existence and uniqueness of a solution of this problem.  In the first two
references above, concrete, physiologically meaningful values for the
parameters are determined and numerical results are presented which are,
however, contradictory.  We used for our experiments four different
parameter sets proposed by McElwain \cite{me:oxygen} and which can be found
in Table~\ref{tab:oxygen}.

\begin{figure}[ht]
  \begin{minipage}{0.5\linewidth}
    \centering
    \includegraphics[width=6cm]{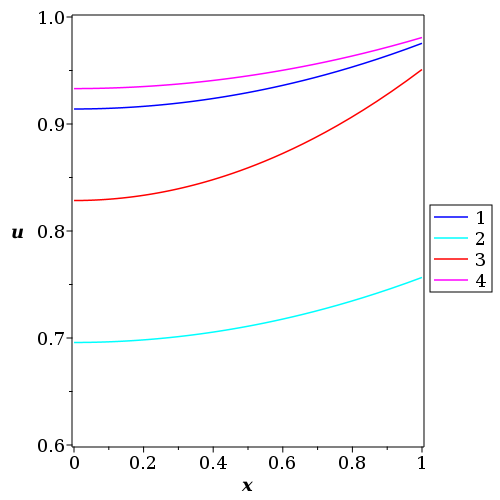}
  \end{minipage}
  \begin{minipage}{0.5\linewidth}
    \centering
    \includegraphics[width=6cm]{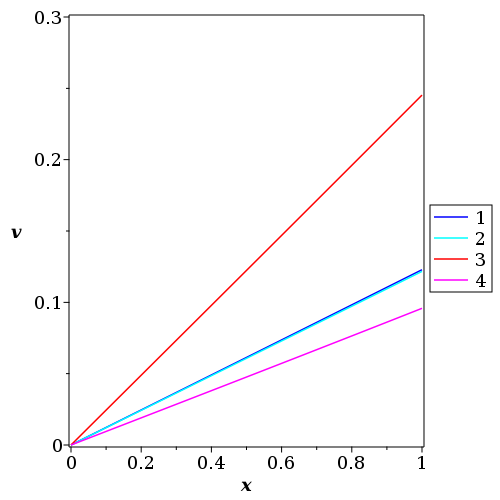}
  \end{minipage}
  \caption{Numerical solutions of the boundary value problem for the oxygen
    uptake model for four different sets of parameters given in
    Table~\ref{tab:oxygen}. Left: oxygen concentration $u(x)$. Right: rate
    of change of oxygen concentration $u'(x)$.}
  \label{fig:oxygen}
\end{figure}

\begin{wraptable}{r}{0.45\linewidth}
  \centering
  \begin{tabular}{c|rrr}
    \hline
    & \multicolumn{1}{c}{$a$} & \multicolumn{1}{c}{$K$} &
    \multicolumn{1}{c}{$\alpha$} \\
    \hline
    1 & $0.38065$ & $0.03119$ & $5$\\
    2 & $0.38065$ & $0.03119$ & $0.5$\\
    3 & $0.76129$ & $0.03119$ & $5$\\
    4 & $0.38065$ & $0.31187$ & $5$\\
    \hline
  \end{tabular}
  \caption{Parameter values for the oxygen uptake model following McElwain
    \cite{hl:oxygen}.}
  \label{tab:oxygen}
\end{wraptable}
In particular for the third parameter set, several authors performed
similar computations starting with Hiltmann und Lory \cite{hl:oxygen}.
Khuri and Sayfy \cite[Ex.~3]{ks:sbvpphys} combined a decomposition method
in the vicinity of the singularity with a collocation method in the rest of
the integration interval.  They provided -- like Hiltmann and Lory --
approximations of $u(x_{i})$ for $x_{i}=i/10$ with $i=0,\dots,10$
\cite[Tbl.~5]{ks:sbvpphys} and compared with results of \c{C}a\u{g}lar et
al.~\cite{cco:bspline}.  It turned out that for the first six digits all
three approaches and our method yield exactly the same result -- a quite
remarkable agreement.  Fig.~\ref{fig:oxygen} provides plots of the oxygen
concentration $u(x)$ and of its rate of change $v(x)=u'(x)$ for all four
different sets of parameters as obtained by our method.  The concentration
plot agrees well with the one given by McElwain \cite[Fig.~1]{me:oxygen}
(and confirmed by Hiltmann und Lory \cite{hl:oxygen}).

Hiltmann and Lory \cite{hl:oxygen} report that they used a sophisticated
implementation of a multiple shooting procedure based on four different
integrators for initial value problems together with a special treatment of
the singularity using both a technique of de Hoog and Weiss
\cite{hw:diffbvpsing} and a Taylor series method (no further details are
given).  They prescribed a tolerance of $10^{-8}$ for their Newton solver
and $10^{-10}$ for the integrator.  By contrast, we used a simple shooting
method with the \textsc{Maple} built-in Runge--Kutta--Fehlberg integrator
and a hand-coded Steffensen method for the nonlinear system with a
tolerance of $10^{-7}$.  This comparison again demonstrates how much
simplicity and robustness one gains by using the associated vector field
for the numerical integration in singular situations.

\subsubsection{Lane--Emden Systems}
\label{sec:les}

Our approach works for systems in the same manner as for scalar equations,
as one still finds a one-dimensional unstable manifold corresponding to
prolonged graph of the solution.  Thus we restrict to just one example of
dimension $d=3$.  We now have to integrate the system \eqref{eq:sysvdlecs}
of dimension $n=2d+1=7$ for the above given initial data.  We used a Newton
method for solving the nonlinear system arising in the shooting method.
Since we had to determine $d=3$ initial conditions via shooting, we had to
augment \eqref{eq:sysvdlecs} by a linear matrix differential equation with
variable coefficients of dimension $7\times 3$.

Campesi et al. \cite{mac:cat} proposed a system of coupled Lane--Emden
equations as model for the combustion of ethanol and ethyl acetate over an
\textsf{MnCu} catalyst using a Langmuir--Hinshelwood--Hougen--Watson
kinetics.  In dimensionless form, the system is given by (see
\cite{mpr:mccp})
\begin{equation}\label{eq:catsys}
  \begin{aligned}
    u''+\frac{2}{x}u' &= \frac{\mu_{u}u}{1+\lambda_{u}u+\lambda_{v}v+\lambda_{w}w}\;,\\
    v''+\frac{2}{x}v' &= \frac{\mu_{v}v-\mu_{u}u}{1+\lambda_{u}u+\lambda_{v}v+\lambda_{w}w}\;,\\
    w''+\frac{2}{x}w' &= \frac{\mu_{w}w}{1+\lambda_{u}u+\lambda_{v}v+\lambda_{w}w}\;,
  \end{aligned}
\end{equation}
where $u$, $v$, $w$ represent (dimensionless) molar concentrations of
ethanol, acetaldehyde and ethyl acetate, respectively.  The boundary
conditions require that at $x=0$ all first derivatives vanish and that at
$x=1$ all concentrations are $1$.

\begin{figure}[ht]
  \begin{minipage}{0.5\linewidth}
    \centering
    \includegraphics[width=6cm]{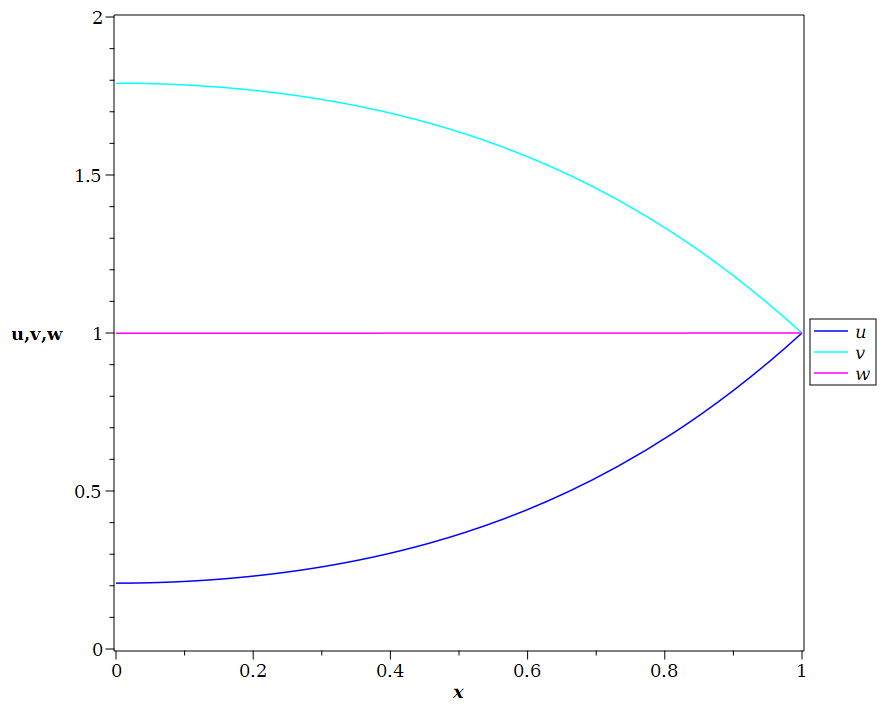}
  \end{minipage}
  \begin{minipage}{0.5\linewidth}
    \centering
    \includegraphics[width=6cm]{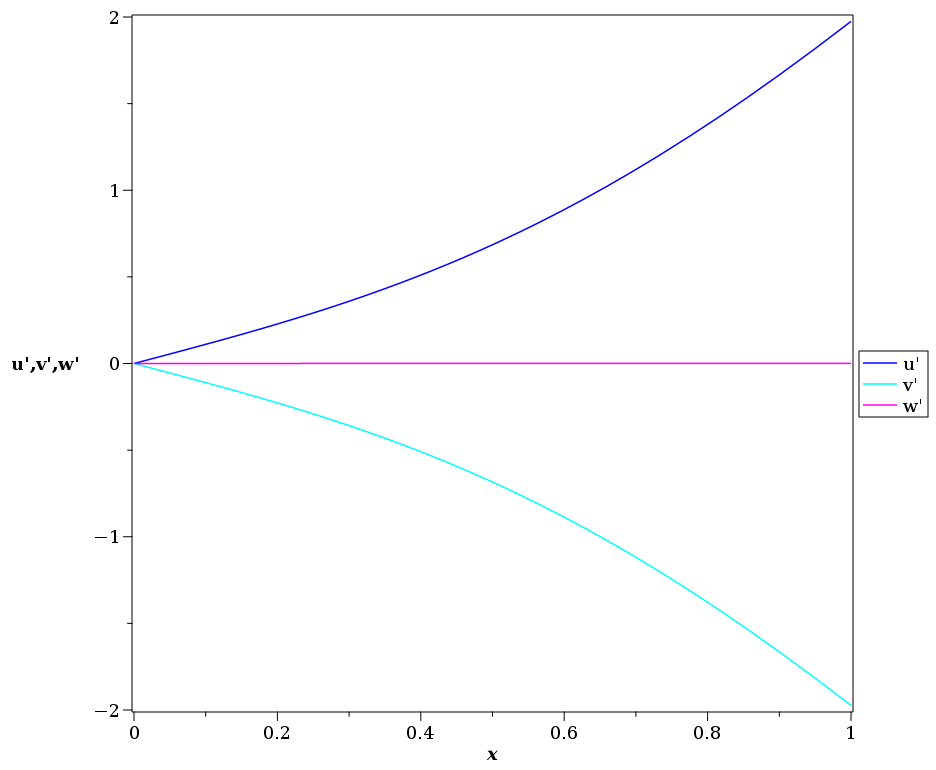}
  \end{minipage}
  \caption{Numerical solutions of the boundary value problem for the
    dimensionless model of the \textsf{\footnotesize MnCu} catalyst. Left:
    concentrations of ethanol, acetaldehyde and ethyl acetate,
    respectively. Right: corresponding rates of change.}
  \label{fig:combust}
\end{figure}

The authors of \cite{mac:cat} used for numerically integrating
\eqref{eq:catsys} an approach developed by essentially the same group
\cite{spb:crn} based on an integral formulation and an $h$-adaptive mesh
procedure.  Unfortunately, \cite{mac:cat} does not provide all the
parameters used in the computations so that it is not possible to compare
with their results.  We used instead for our experiments data given in
\cite{mpr:mccp} (employing a modified Adomian decomposition method).
However, the plots given there are not correct, as apparently wrong
differential equations were used -- at least in the \textsc{Matlab} code
presented in the appendix.  We compared with analogous \textsc{Matlab}
computations using the right differential equations and again
\texttt{pdepe} as a numerical solver and obtained an excellent agreement.
Figure~\ref{fig:combust} presents solution curves for the values
$\mu_{u}=30$, $\mu_{v/w}=0.01$, $\lambda_{u}=3$ and $\lambda_{v/w}=0.1$
used in \cite{mpr:mccp}.

\section{Thomas--Fermi Equation}
\label{sec:tf}

\subsection{Majorana Transformation}
\label{sec:tfmaj}

The Thomas--Fermi equation \eqref{eq:tf} belongs also to the class
\eqref{eq:sql2e}, but with
\begin{equation}
  g(x)=\sqrt{x}\;,\qquad f(x,u,u')=\sqrt{u^{3}}\;.
\end{equation}
The initial condition $u(0)=1$ leads to a rather different situation as for
the Lane--Emden equation: the implicit form of the Thomas--Fermi equation
entails that the only points on $\mathcal{R}_{2}$ which project on $x=0$
are of the form $\rho=(0,0,u_{1},u_{2})$ with arbitrary values $u_{1}$,
$u_{2}$.  Hence no solution satisfying the above initial condition can be
twice differentiable at $x=0$.  Solutions with a higher regularity exist
only for the initial condition $u(0)=0$ which has no physical relevance.

Any point of the form $\rho=(0,1,u_{1})$ is an \emph{improper} impasse
point.  The vector field $Y$ defined by \eqref{eq:Yvd} does not vanish at
such points but takes the form $\partial_{u'}$ and it is not Lipschitz
continuous there.  While Peano's theorem still asserts the existence of
solutions, we cannot apply the Picard--Lindel\"of theorem to guarantee
uniqueness.  We could rescale $Y$ by some function like $x$ which does not
change its trajectories for obtaining an everywhere differentiable vector
field $\tilde{Y}=xY$.  Now all points of the above form are stationary
points.  But the Jacobian of $\tilde{Y}$ has $0$ as a triple eigenvalue at
them making it hard to analyse the local phase portrait.

We use therefore a different approach.  As Esposito \cite{se:mtf} reported
only in 2002, Majorana proposed already in 1928 a differential
transformation to a new independent variable $t$ and a new dependent
variable $v$ of the form
\begin{equation}\label{eq:majtrafo}
  t=144^{-1/6}x^{1/2}u^{1/6}\;,\qquad
  v=-(16/3)^{1/3}u^{-4/3}u'\;.
\end{equation}
This at first sight rather miraculous transformation stems from a
particular kind of scaling symmetry \cite{se:mtde}.  A computation detailed
in \cite{se:mtf} shows that if it is applied to any solution of the
Thomas--Fermi equation \eqref{eq:tf}, then the transformed variables
satisfy the reduced equation
\begin{equation}\label{eq:tfmaj}
  (1-t^{2}v)\frac{dv}{dt}=8(tv^{2}-1)\;.
\end{equation}
The boundary condition \eqref{eq:tfbcinf}, i.\,e.\
$\lim_{x\to\infty}u(x)=0$, translates into the condition
$v(1)=1$.\footnote{The Majorana transformation is \emph{not} bijective.  A
  well-known solution of the Thomas--Fermi equation already given by Thomas
  \cite{lht:atom} is $u_{s}(x)=144x^{-3}$.  It does not satisfy the left
  boundary condition, as it is not even defined for $x=0$, but the
  asymptotic condition at infinity.  One easily verifies that any point of
  the form $\bigl(x,u_{s}(x),u'_{s}(x)\bigr)$ is mapped into the point
  $(1,1)$.}  We will see below that the thus defined singular initial value
problem for \eqref{eq:tfmaj} possesses \emph{two} solutions.  Only one of
them is also defined for $t=0$ and thus is the physically relevant one.  It
follows from \eqref{eq:majtrafo} that the initial slope $u'(0)$ for the
Thomas--Fermi equation is obtained from a solution of \eqref{eq:tfmaj} by
\begin{equation}\label{eq:tfslope}
  u'(0)=-(3/16)^{1/3}v(0)\;.
\end{equation}

The reduced equation \eqref{eq:tfmaj} is quasi-linear and of first order.
Opposed to the Lane--Emden equations, it is not semi-linear.  Thus singular
behaviour does not simply occur at specific $t$-values.  Instead it appears
whenever a solution graph contains a point $(t,v)$ with $t^{2}v=1$.
Nevertheless, one can apply the same kind of approach.  One first computes
a vector field $X$ living on the hypersurface
$\mathcal{R}_{1}\subset\mathcal{J}_{1}$ defined by \eqref{eq:tfmaj} and
spanning there the Vessiot distribution.  Then one projects $X$ to the jet
bundle $\mathcal{J}_{0}$ and obtains there the vector field
\begin{equation}\label{eq:tfYred}
  Y_{\mathrm{red}}=(t^{2}v-1)\partial_{t}+8(1-tv^{2})\partial_{v}\;.
\end{equation}
As we are now on $\mathcal{J}_{0}$, one-dimensional invariant manifolds of
$Y_{\mathrm{red}}$ which are transversal can be directly identified with
the graphs of solutions of \eqref{eq:tfmaj}.  Our initial point $(1,1)$ is
a proper impasse point where $Y_{\mathrm{red}}$ vanishes.

\begin{wrapfigure}{r}{0.45\linewidth}
  \centering
  \includegraphics[width=6cm]{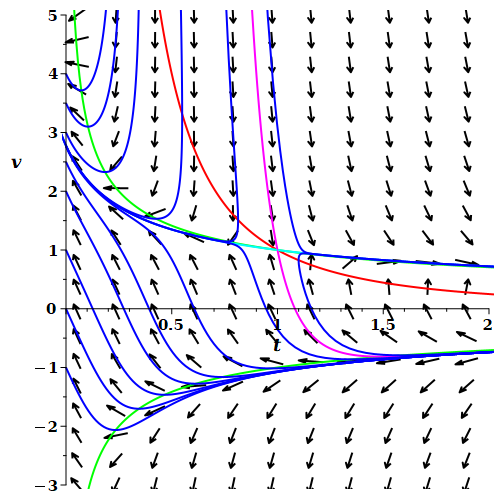}
  \caption{Phase portrait of the vector field associated to the reduced
    system~\eqref{eq:tfmaj}. The unstable manifold is shown in cyan, the
    stable manifold in magenta.}
  \label{fig:ppmaj}
\end{wrapfigure}

Fig.~\ref{fig:ppmaj} shows the phase portrait of the vector field
$Y_{\mathrm{red}}$.  It has $(1,1)$ as its only stationary point.  The plot
shows in blue some integral curves.  Most, but not all of them can be
considered as the graphs of solutions of \eqref{eq:tfmaj}.  The plot also
contains in red the $t$-nullcline given by $v=1/t^{2}$ -- which is
simultaneously the singular locus of \eqref{eq:tfmaj} -- and in green the
$v$-nullcline given by $v=\pm1/\sqrt{t}$.  The integral curves that cross
the $t$-nullcline show at the intersection a turning point behaviour, as
the $t$-component of $Y_{\mathrm{red}}$ changes its sign there.  If
$(t_{i},v_{i})$ is such an intersection point, then it splits the
corresponding integral curve into \emph{two} solution graphs where both
solutions are defined only for $t<t_{i}$, as they both loose their
differentiability at $t=t_{i}$.  With traditional numerical methods applied
to \eqref{eq:tfmaj}, it would be difficult to determine these solutions; as
integral curves of $Y_{\mathrm{red}}$ they are trivial to obtain
numerically.

The Jacobian of $Y_{\mathrm{red}}$ at the stationary point $(1,1)$ is the
matrix
$J=\left(\begin{smallmatrix} 2 & 1 \\ -8 & -16 \end{smallmatrix}\right)$
with eigenvalues $-7\pm\sqrt{73}\approx(1.544,-15.544)$.  Thus we are
dealing with a saddle point.  The unstable and the stable manifold shown in
Fig.~\ref{fig:ppmaj} in cyan and magenta, resp., correspond to the above
mentioned two solutions of the initial value problem with $v(1)=1$. There
cannot exist any additional solutions, as there are no further invariant
manifolds entering or leaving the saddle point.  One sees that in the
positive quadrant the stable manifold cannot cross the nullclines outside
of the saddle point and hence can never reach the $v$-axis.

Thus we may conclude that the part of the unstable manifold between the
$v$-axis and the stationary point corresponds to the unique solution
$u_{\infty}$ of the boundary value problem with the condition
\eqref{eq:tfbcinf}.  The abscissa of the intersection of the unstable
manifold with the $v$-axis determines via \eqref{eq:tfslope} the critical
initial slope $\omega$ (see below for numerical values).  The existence of
such a unique solution for this specific boundary value problem was proven
in 1929 by Mambriani \cite{am:edo2} (see also the discussion in
\cite{eh:tf2}).

It will turn out that the integral curves to the right of the stable
manifold have no relevance for our problem.  The integral curves to the
left of it and above the unstable manifold correspond to solutions of the
boundary value problem with the condition \eqref{eq:tfbcion}, i.\,e.\
solutions with a zero, whereas the integral curves below the stable
manifold lead to solutions for \eqref{eq:tfbccryst}.  This can be deduced
from their intersections with the $v$-axis and \eqref{eq:tfslope}.

Much of the literature on numerically solving the Thomas--Fermi equation is
concerned with the solution $u_{\infty}$ of \eqref{eq:tfbcinf} defined on
the semi-infinite interval $[0,\infty)$ and concentrates on the
determination of the critical slope~$\omega$.  Most solutions reported in
the literature are either shown only on rather small intervals $[0,x_{0}]$
with typically $x_{0}<10$ or clearly deteriorate for larger $x$.  One
reason for this effect is surely that many approaches are based on some
kind of series expansion.  Another, more intrinsic reason becomes apparent
from the phase portrait in Figure~\ref{fig:ppmaj}.  As the sought solution
corresponds to a branch of the unstable manifold of the saddle point
$(1,1)$, even small errors close to the saddle point (corresponding to
points with large $x$ coordinates) are amplified by the dynamics and the
numerical solutions tend to diverge from a finite limit.

By contrast, our approach to determine $u_{\infty}$ leads to the standard
problem of determining a branch of the unstable manifold of a stationary
point -- a task which can be performed numerically very robustly and
efficiently.  As the positive eigenvalue has about the tenfold magnitude of
the negative one, trajectories approach the unstable manifold very fast
which ensures a high accuracy.

Following Majorana, Esposito \cite{se:mtf} (and subsequent authors)
determines a series solution of the initial value problem $v(1)=1$ for
\eqref{eq:tfmaj}.  In the first step, one obtains a quadratic equation with
two solutions.  Esposito then argues that one should take the smaller
solution, as this was a perturbation calculation which is not a convincing
argument.  The reduced initial value problem has \emph{two} solutions.  As
one can see in Figure~\ref{fig:ppmaj}, the second solution corresponding to
the stable manifold grows very rapidly.  Therefore it is not surprising
that several authors suspected that the second solution of the quadratic
equation leads to a divergent power series and thus could be discarded.
However, a second solution to the initial value problem does exist,
although it seems that it cannot be determined with a power series ansatz.
But as already discussed above, $u_{\infty}$ is nevertheless unique and
corresponds to the unstable manifold.

For the series solution, one expands around $t=1$ and makes the ansatz
$v(t)=\sum_{i=0}^{\infty}a_{i}(1-t)^{i}$.  The initial condition yields
$a_{0}=1$ and for the arising quadratic equation for $a_{1}$ one chooses
the root\footnote{This value is related to the spectrum of the Jacobian of
  the vector field $Y_{\mathrm{red}}$: $-a_{1}$ is the slope of the tangent
  space of the unstable manifold at the saddle point.  This is not
  surprising, as the tangent space is the linear approximation of the
  solution.}  $a_{1}=9-\sqrt{73}\approx 0.456$.  After lengthy computations
sketched in \cite{se:mtf}, one obtains the following recursive expression
for the remaining coefficients with $i>1$:
\begin{multline}\label{eq:majseries}
  a_{i}=\frac{1}{2(i+8)-(i-1)a_{1}}
  \Biggl[(i+6)a_{1}a_{i-2}+\Bigl((i+7)-2(i+3)a_{1}\Bigr)a_{i-1}+{}\\
  \sum_{j=1}^{i-2}\Bigl((j+1)a_{j+1}-2(j+4)a_{j}+(j+7)a_{j-1}\Bigr)a_{i-j}\Biggr]\;.
\end{multline}
Setting $t=0$ yields for the critical slope the series representation
\begin{equation}\label{eq:slopeseries}
  \omega=-\left(\frac{3}{16}\right)^{1/3}\sum_{i=0}^{\infty}a_{i}\;,
\end{equation}
which can be evaluated to arbitrary precision.

To obtain whole solutions $u(x)$, one must be able to transform back from
the variables $(t,v)$ to the original variables $(x,u)$.  Esposito
\cite{se:mtf} exhibited a convenient method for this.  We express the
solution in parametric form using $t$ as parameter: $x=x(t)$ and $u=u(t)$.
Then we make the ansatz
\begin{equation}\label{eq:ansuw}
  u(t)=\exp{\left\{\int_{0}^{t}w(\tau)d\tau\right\}}
\end{equation}
with $w$ a yet to be determined function.  Assuming $x(t=0)=0$, this ansatz
automatically satisfies the initial condition $u(x=0)=1$.  Using the
transformation \eqref{eq:majtrafo}, one can show that
$w(t)=\tfrac{6tv(t)}{t^{2}v(t)-1}$ and that $x(t)$ can be expressed via
$w(t)$ as
\begin{equation}\label{eq:xt}
  x(t)=144^{1/3}t^{2}\exp{\left\{-\frac{1}{3}\int_{0}^{t}w(\tau)d\tau\right\}}
\end{equation}
(which shows that indeed $x(0)=0$).  Esposito \cite{se:mtf} proposed to
enter the above determined series solution for $v(t)$ into these formulae
and to compute this way a series expansion of $u_{\infty}$.  This requires
essentially one quadrature.

\subsection{Numerical Results}
\label{sec:tfnum}

We refrain from citing the many papers written on computing $u_{\infty}$
and in particular $\omega$ and instead refer only to \cite{pd:tfrcf,fg:mtf}
both listing a large number of approaches with references.  We emphasise
again that our main point is to show that the geometric theory allows us --
here in combination with the Majorana transformation -- to translate a
singular problem into basic tasks from the theory of dynamical systems
which can be easily solved by standard methods.

\subsubsection{The ``Critical'' Solution $u_{\infty}$ and the Critical
  Slope $\omega$}

We consider first the problem of only determining the initial slope
$\omega$ belonging to the solution $u_{\infty}$ for \eqref{eq:tfbcinf}.
With classical approaches, this is a non-trivial problem and in the
literature one often finds values with a very low number of correct digits.
Using our geometric approach, we can determine $\omega$ to (almost) any
desired precision in about 10 lines of \textsc{Maple} code.  We write the
dynamical system corresponding to the vector field $Y_{\mathrm{red}}$
defined by \eqref{eq:tfYred} as
\begin{equation}\label{eq:dynsysYred}
  \frac{dt}{ds}=t^{2}v-1\;,\qquad \frac{dv}{ds}=8(1-tv^{2})\;,
\end{equation}
i.\,e.\ we determine integral curves of $Y_{\mathrm{red}}$ in parametric
form $\bigl(t(s),v(s)\bigr)$.  As discussed above, the sought trajectory
corresponds to the unstable manifold of the saddle point~$(1,1)$.  An
eigenvector for the positive eigenvalue $\lambda=-7+\sqrt{73}$ is given by
$\mathbf{e}=\bigl(1,\ -9+\sqrt{73}\bigr)^{T}$ and we denote by
$\hat{\mathbf{e}}=(e_{1},\ e_{2})^{T}$ the corresponding normalised vector.
Then we choose as initial point for a numerical integration
$t(0)=1+\epsilon e_{1}$ and $v(0)=1+\epsilon e_{2}$ with $\epsilon>0$ some
small number (we typically used $10^{-3}$ or $10^{-4}$, but this had no
effect on the obtained slope) and integrated until $t(s)=0$ for $s=s_{0}$.
Finally, we obtain $\omega$ from $v(s_{0})$ via \eqref{eq:tfslope}.

We control the precision with an integer parameter $N$ specifying that the
numerical integration of \eqref{eq:dynsysYred} should take place with an
absolute and relative error of $10^{-N}$ and that for this purpose
\textsc{Maple} should compute with $N+5$ digits.  In a recent work,
Fern\'andez and Garcia \cite{fg:mtf} determined $\omega$ based on the first
5000 terms of the Majorana series \eqref{eq:slopeseries} to a precision of
several hundred digits.  This is by far the best approximation available
and our reference solution.

\begin{wraptable}{r}{0.45\linewidth}
  \centering
  \begin{tabular}{rrr}
    \hline
    \multicolumn{1}{c}{tolerance} & \multicolumn{1}{c}{rel. error} &
    \multicolumn{1}{c}{time} \\
    \hline
    $10^{-5\phantom{0}}$ & $3.2\cdot 10^{-6\phantom{0}}$ & $0.6$\\
    $10^{-10}$ & $7.3\cdot 10^{-12}$ & $0.6$\\
    $10^{-15}$ & $5.5\cdot 10^{-17}$ & $2.7$\\
    $10^{-20}$ & $5.3\cdot 10^{-22}$ & $22.7$\\
    $10^{-25}$ & $5.5\cdot 10^{-27}$ & $231.5$\\
    \hline
  \end{tabular}
  \caption{Relative error and computation time in seconds for different tolerances.}
  \label{tab:ft}
\end{wraptable}
Our numerical results are summarised in Table~\ref{tab:ft}.  Our relative
error is always smaller than the prescribed tolerance.  For smaller
tolerances, the computational effort is rapidly increasing and on a laptop
we needed for 25 digits less than 4 minutes.  We made no effort to optimise
the computations.  For example, we are using the default integration method
of \textsc{Maple} (a Runge--Kutta--Fehlberg method of order 4/5 with a
degree four interpolant), although a higher order scheme would probably be
more efficient (\textsc{Maple} offers such schemes -- but not in
combination with the automated root finding used in our code).
Nevertheless, one may conclude that for practically relevant precisions,
our geometric approach combined with the Majorana transformation provides
very accurate results fast and almost effortless.

\begin{table}[ht]
  \centering
  \begin{tabular}{cccccc}
    \hline
    terms & 10 & 20 & 30 & 40 & 50 \\ 
    rel. err. & $5.8\cdot 10^{-2}$ & $6.7\cdot 10^{-3} $&
    $8.3\cdot 10^{-4}$ & $1.1\cdot 10^{-4}$ & $1.4\cdot 10^{-5}$ \\
    \hline
    terms & 60 & 70 & 80 & 90 & 100\\
    rel. err. & $1.9\cdot 10^{-6}$ &  $2.7\cdot 10^{-7}$
    & $3.7\cdot 10^{-8}$ & $4.4\cdot 10^{-9}$ & $0$\\
    \hline
  \end{tabular}
  \caption{Relative error for different truncation degrees of the Majorana series.}
  \label{tab:ftm}
\end{table}

Fern\'andez and Garcia \cite{fg:mtf} analyse also the convergence rate of
the Majorana series \eqref{eq:slopeseries} and consider it as fast (see
also the comments by Esposito \cite{se:mtf}).  We compared for a relative
small accuracy, \textsc{Maple} hardware floats with 10 digits, the value
for the initial slope obtained with our approach with the approximations
delivered by various truncations of the series.  Somewhat surprisingly, our
approach gets all 10 digits right, despite the considerably higher
tolerances ($10^{-6}$) used by the integrator.  Table~\ref{tab:ftm}
contains the approximations obtained by evaluating the first $N$ terms of
the Majorana series \eqref{eq:slopeseries}.  One needs 100 terms for a
similarly accurate result.  On average, one needs 10 more terms for one
additional digit corresponding to a linear convergence as already
theoretically predicted in \cite{se:mtf,fg:mtf}.  This observation also
roughly agrees with the fact that Fern\'andez and Garcia used 5000 terms
for obtaining about 500 digits \cite{fg:mtf}.

For determining the whole solution $u_{\infty}(x)$ instead of only the
critical slope $\omega=u_{\infty}'(0)$, we have to perform a transformation
back from the variables $(t,v)$ to $(x,u)$.  We described above Esposito's
approach for this.  For a purely numerical computation instead of series
expansions, we modify it in a way which fits nicely into our approach.  We
introduce as Esposito \cite{se:mtf} the function
\begin{equation}\label{eq:defI}
  \mathcal{I}(t)=\int_{0}^{t}\frac{\tau v(\tau)}{1-\tau^{2}v(\tau)}d\tau\;.
\end{equation}
We then express $\mathcal{I}(t)$ as a function of the parameter $s$ which
we use to parametrise solution curves.  If $s_{0}$ is the (first) parameter
value satisfying $t(s_{0})=0$, then an elementary application of the
substitution rule yields
\begin{equation}\label{eq:Is}
  \mathcal{I}(s)=-\int_{s_{0}}^{s}t(\sigma)v(\sigma)d\sigma\;,
\end{equation}
which immediately implies that $\mathcal{I}$ satisfies the differential
equation $\tfrac{d\mathcal{I}}{ds}=-tv$ by which we augment the system
\eqref{eq:dynsysYred}.  We thus obtain a \emph{free} boundary value problem
for the augmented system, as the function $\mathcal{I}(s)$ satisfies the
condition $\mathcal{I}(s_{0})=0$ with the a priori unknown value $s_{0}$.
As usual, we consider $s_{0}$ as an additional unknown function and
introduce the rescaled independent variable $\sigma=s/s_{0}$.  Then we
finally obtain the following two-point boundary value problem with
non-separated boundary conditions
\begin{equation}\label{eq:tfexsys}
  \begin{aligned}
    \frac{dt}{d\sigma} &= s_{0}(t^{2}v-1)\;, & t(0) &=1+\epsilon e_{1}\;, &
    t(1) &= 0\;,\\[0.3\baselineskip]
    \frac{dv}{d\sigma} &= 8s_{0}(1-tv^{2})\;, & v(0) &= 1+\epsilon
    e_{2}\\[0.3\baselineskip]
    \frac{d\mathcal{I}}{d\sigma} &= -s_{0}tv\;, &&& \mathcal{I}(1)
    &=0\;,\\[0.3\baselineskip]
    \frac{ds_{0}}{d\sigma} &= 0\;.
  \end{aligned}
\end{equation}
Once this boundary value problem is solved, \eqref{eq:ansuw} and
\eqref{eq:xt} imply that parametrisations of the graph of $u_{\infty}(x)$
are given by
\begin{equation}\label{eq:tfparasol}
  x(\sigma) =
  144^{1/3}t(\sigma)^{2}\exp{\bigl(2\mathcal{I}(\sigma)\bigr)}\;,\quad
  u(\sigma)= \exp{\bigl(-6\mathcal{I}(\sigma)\bigr)}\;.
\end{equation}

\begin{wrapfigure}{r}{0.45\linewidth}
  \centering
  \includegraphics[width=6cm]{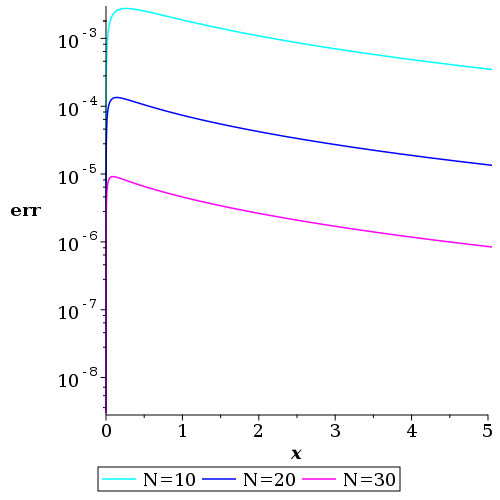}
  \caption{Comparison of values obtained via \eqref{eq:tfparasol} and
    Majorana's series for different numbers $N$ of terms.\label{fig:tfsolerr}}
\end{wrapfigure}
We implemented this approach in \textsc{Maple} using the built-in solver
for boundary value problems which could handle \eqref{eq:tfexsys} without
problems.  We compared the results with solutions obtained via Majorana's
series, i.\,e.\ following Esposito \cite{se:mtf}, we entered a given number
$N$ of terms into the integral defining $\mathcal{I}$ and performed a
numerical integration.  Fig.~\ref{fig:tfsolerr} shows on a logarithmic
scale the absolute difference between our curve
$\bigl(x(\sigma),u(\sigma)\bigr)$ and the curves computed via the series
for different values of $N$.  Obviously, our results are in an excellent
agreement with the series solutions.  The fact that all error curves have
their maximum close to $x=0$ is easy to explain.  As the expansion point of
the series corresponds to $x=\infty$ (i.\,e.\ $t=1$), the series solutions
become less accurate the closer one gets to $x=0$; at $x=0$ of course no
error occurs, as this value is fixed by an initial condition.  We did not
make an extensive comparison of computation times.  But plotting the series
solution for $N=10$ over the interval $[0,5]$ required more than 10 times
as much computation time than solving above boundary value problem
demonstrating again the efficiency of our approach.

We mentioned already above that in the literature results are often
presented only for rather small values of $x$, although the solution is
defined for all non-negative real numbers.  One exception is Amore et
al.~\cite[Tbl.~3/4]{abf:acctf} who used a Pad\'e--Hankel method and
asymptotic expansions to present highly accurate values of the solution
$u_{\infty}(x)$ and its first derivative $u_{\infty}'(x)$ up to $x=400$.

\begin{table}[ht]
  \centering
  \begin{tabular}{rrr}
    \hline
    \multicolumn{1}{c}{$x$} &  \multicolumn{1}{c}{$u_{\infty}(x)$} &
   \multicolumn{1}{c}{$u_{\infty}'(x)$}  \\ 
    \hline
    $0$ & $1$ & $-1.58807101687867$ \\
    $10$ & $0.0243142929534589$ & $-0.00460288186903816$ \\
    $50$ & $0.000632254782228818$ & $-0.0000324989019998445$ \\
    $100$ & $0.000100242568239745$ & $-2.73935106365787 \cdot 10^{-6}$ \\
    $150$ & $0.0000326339644201454$ & $-6.09139947257267 \cdot 10^{-7}$ \\
    $200$ & $0.0000145018034835377$ & $-2.05753231409599 \cdot 10^{-7}$ \\
    $250$ & $7.67729076668264 \cdot 10^{-6}$ & $-8.78946798702223 \cdot 10^{-8}$ \\
    $300$ & $4.54857195240339 \cdot 10^{-6}$ & $-4.36594961733055 \cdot 10^{-8}$ \\
    $350$ & $2.91510210708972 \cdot 10^{-6}$ & $-2.40920109677041 \cdot 10^{-8}$ \\
    $400$ & $1.97973262954641 \cdot 10^{-6}$ & $-1.43668230750324 \cdot 10^{-8}$ \\
    \hline
  \end{tabular}
  \caption{Solution values $u_{\infty}(x)$ and derivative values
    $u_{\infty}'(x)$ for large $x$.}\label{tab:ftlx}
\end{table}

Table~\ref{tab:ftlx} contains similar values obtained with our approach.
For determining the values of $u_{\infty}'(x)$, we must augment
\eqref{eq:tfparasol} by an equation for $u'(\sigma)$, i.\,e.\ we must
extend the parametrisation to the prolonged graph.  By a straightforward
application of the chain rule, one obtains
\begin{equation}\label{eq:ftparader}
  u'(\sigma)=-3\cdot 144^{-1/3}v(\sigma)\exp{\bigl(-8\mathcal{I}(\sigma)\bigr)}\;.
\end{equation}
To compile such a table, one must then determine for each $x$ the
corresponding value of the parameter $\sigma$ via the solution of a
nonlinear equation.  Nevertheless, the complete computation of the values
at the ten points contained in the table required only about 0.1 seconds.
Amore et al.  \cite{abf:acctf} claim that in their tables all digits are
correct.  Assuming that this is indeed the case, we can conclude that we
obtained with minimal effort for each value of $x$ at least eight correct
digits for $u_{\infty}(x)$ and seven correct digits for $u_{\infty}'(x)$.
Given the settings for the tolerances of our integrator and the use of
hardware floats with only 10 digits, these results demonstrate again a very
remarkable precision and efficiency of our approach.  As large values of
$x$ correspond to small values of $\sigma$ and thus to values of $t$ close
to $1$, one may have to choose a smaller value of $\epsilon$ for very large
values of $x$.  The largest value appearing in above table, $x=400$,
corresponds to $\sigma\approx 0.35$ and $t\approx 0.9789$.  We chose for
our numerical calculation the value $\epsilon=10^{-3}$ and thus used as
right end of the approximated unstable manifold instead of the saddle point
$(1,1)$ the point $(t_{1},v_{1})\approx (0.9978,1.001)$.  For $x=400$, one
may say that we are still sufficiently far away from this point, but for
larger values of $x$ one should probably start working with a smaller value
of $\epsilon$ which will increase the computation time, as the dynamics is
very slow so close to a stationary point.

\subsubsection{Other Solutions}

So far, we only considered the particular solution $u_{\infty}$ (which has
attracted the most attention in the literature).  In Fig.~\ref{fig:ppmaj}
we presented the phase portrait for the Majorana transformed Thomas--Fermi
equation.  Using a slight modification (and simplification) of the above
described backtransformation via the solution of an extended differential
system, we can also obtain a ``phase portrait'' of the original
Thomas--Fermi equation, i.\,e.\ we compute solutions for different values
of the initial slope $u'(0)$ keeping the initial condition $u(0)=1$.  While
the Majorana transformation itself is valid for any solution of the
Thomas--Fermi equation, our ansatz for the back transformation has encoded
this second initial condition (one could easily adapt to a different value
$u(0)=c$ by multiplying \eqref{eq:ansuw} with the constant $c$).  According
to \eqref{eq:tfslope}, each value of $u'(0)$ corresponds uniquely to a
value of $v(0)$.  We now take the vector field $-Y_{\mathrm{red}}$ and use
a parametrisation such that $s=0$ corresponds to $t=0$ (and thus also
$x=0$).  This leads to the following augmented initial value problem:
\begin{equation}\label{eq:tfppsys}
  \begin{gathered}
    \frac{dt}{ds} = 1-t^{2}v\;,\quad  \frac{dv}{ds} = 8(tv^{2}-1)\;,\quad
    \frac{d\mathcal{I}}{ds} = tv\;,\\[0.3\baselineskip]
    t(0) = 0\;,\quad v(0) = v_{0}\;,\quad \mathcal{I}(0) = 0\;.
  \end{gathered}
\end{equation}
Its solutions are then transformed into $x$- and $u$-coordinates via
\eqref{eq:tfparasol}. 

\begin{wrapfigure}{r}{0.45\linewidth}
  \centering
  \includegraphics[width=6cm]{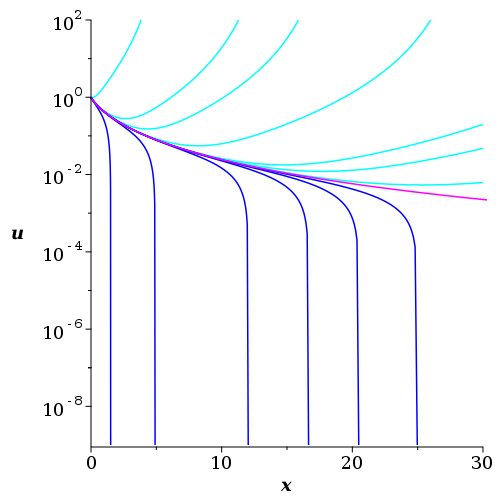}
  \caption{Solutions of the Thomas--Fermi equation with $u(0)=1$ and
    different $u'(0)$ using a logarithmic scale for $u$.  The curve in
    magenta shows $u_{\infty}$.}
  \label{fig:pptf}
\end{wrapfigure}

Fig.~\ref{fig:pptf} shows that the solution $u_{\infty}$ vanishing at
infinity acts as a kind of ``separatrix''.  The solutions above it, i.\,e.\
with an initial slope higher than $\omega$, pass through a minimum and then
grow faster than exponentially (note the logarithmic scale).  The solutions
below it approach rapidly zero, reaching it at a finite value of $x$
(recall that the separatrix reaches zero at infinity).  It turns out that
around the critical value $\omega$, the trajectories are rather sensitive
with respect to the initial slope.  For some of the curves shown in
Fig.~\ref{fig:pptf}, $u'(0)$ differs only in the fifth or sixth digit.  For
the curves approaching zero, it is also non-trivial to determine the exact
location of the zero, as here $v$ goes towards infinity. In our
computations, we actually integrated only until some threshold like
$10^{-8}$.  Probably a ``hybrid'' approach using \eqref{eq:tfppsys} only to
get away from the singularity at $x=0$ and applying afterwards a standard
integrator to the Thomas--Fermi equation would be a good alternative.

For solving concrete boundary value problems with boundary conditions of
the form \eqref{eq:tfbccryst} or \eqref{eq:tfbcion}, resp., for given
values of $a$ or $b$, resp., one can use an adapted version of a shooting
method.  Starting with an initial guess $v_{0}$ for the unknown value of
$v(0)$ for the sought solution, one integrates the initial value problem
\eqref{eq:tfppsys} until a condition of the desired form is satisfied.
However, in general, the condition will be satisfied at a wrong position
$a^{*}$ or $b^{*}$, resp. Using a bisection, one modifies $v_{0}$ until one
is sufficiently close to the actually prescribed values. As in both cases,
Fig.~\ref{fig:pptf} shows that there is a monotone relation between $v_{0}$
and $a^{*}$ or $b^{*}$, resp., it is always clear in which direction one
has to change $v_{0}$. But for larger values of $a$ or $b$, one gets again
into areas where very small changes in $v_{0}$ lead to significant changes
in $a^{*}$ or $b^{*}$, resp.  Despite this sensitivity, the approach worked
in tests very well for $a\leq27$ and $b\leq30$.

\section{Conclusions}

The Lane--Emden and the Thomas--Fermi equation are prototypical examples
for ordinary differential equations with singularities.  Their
singularities are determined by a specific value of the independent
variable: $x=0$.  Any initial or boundary value problem with conditions
prescribed at $x=0$ cannot be tackled by standard methods and this concerns
both theoretical and numerical studies.

The Lane--Emden equations fit into the framework of so-called
\emph{Fuchsian equations} (see e.\,g.\ \cite{sk:fuchs}), i.\,e.\ equations
of the form $Lu=f(x,u)$ where $L$ is a linear differential operator of
Fuchsian type and where only the right hand side may contain nonlinear
terms.  For the theoretical treatment of such equations, some form of
quasilinearisation is often fruitful, as it allows to use the far developed
theory of the linear counterpart $Lu=\tilde{f}(x)$.  For example, the
existence and uniqueness proof for boundary value problems for
(generalised) Lane--Emden equations given in \cite{cs:exssbvp} follows such
strategy.  For the numerical integration, \cite{kkw:singivp} presents
methods for first- and second-order systems of this particular form.

A key consequence of this special structure is the above mentioned location
of the singularities depending only on $x$ which facilitates the design of
specialised numerical methods.  Therefore it is not surprising that so many
different techniques have been proposed in the literature.  Our approach is
independent of such a special form, as one can see from our treatment of
the Thomas--Fermi equation based on the reduced equation \eqref{eq:tfmaj}.
The location of its singularities depends on $t$ and $v$ making an
integration with standard numerical methods more difficult.  By contrast,
our approach can handle all forms of quasilinear problems.

In some computations related to the Thomas--Fermi equations, we encountered
problems, for example when computing $u_{\infty}(x)$ for very large values
of $x$ or when $u(x)$ approaches zero.  In the first case, the reason lies
in an often highly nonlinear relationship between the variable $t$ used in
the reduced system and the variable $x$ where ``microscopic'' changes in
$t$ may correspond to huge differences in $x$.  In the second case, $u$ can
approach $0$ only when $v$ tends towards infinity.  In both cases, one
could probably extend the applicability of our method by a rescaling of the
reduced equation.  For computing $u_{\infty}$ for large $x$, an
alternative, semianalytic approach would consist of determining a higher
order approximation of the unstable manifold close to the saddle point
$(1,1)$ -- in fact, the Majorana series is nothing else than such an
approximation.  This could lead to very accurate values even for extremely
large values of $x$.

One may wonder why we used in the case of the Lane--Emden equations the
shooting method for boundary value problems and not also a formulation as
free boundary value problem as for the Thomas--Fermi equation.  In both
cases, one faces the problem that at one boundary one has to deal with a
two-dimensional plane of stationary points and that the boundary conditions
enforce that one end point of the solution trajectory lies on this plane.
In the case of the Thomas--Fermi equation, we resolved this problem by
moving a bit in the direction of the unstable eigenspace.  This was
possible, as this direction is the same for all points on the plane.  In
the case of the Lane--Emden equations, the direction of the unstable
eigenspace depends on the value $u(0)$ and thus differs for different
points on the plane.  Probably one could adapt typical approaches to
boundary value problems like collocation methods to this dependency.  But
as our emphasis in this paper lies on the use of standard methods, we
refrained from studying this possibility in more details.  Furthermore, the
simple shooting method works very well and reliable for this class of
problems.

\section*{Acknowledgements}

This work was performed within the Research Training Group \emph{Biological
  Clocks on Multiple Time Scales} (GRK 2749) at Kassel University funded by
the German Research Foundation (DFG).

\section*{Data}

All computations presented in this article (and a few more) were performed
within \textsc{Maple} worksheets.  These are freely available at Zenodo
using the DOI \url{https://doi.org/10.5281/zenodo.10064926}.

\bibliographystyle{elsarticle-num}
\bibliography{QuasiNum}

\end{document}